\magnification=\magstep1
\input amstex
\documentstyle{amsppt}
\pageheight{47pc}
\pagewidth{33pc}
\NoBlackBoxes
\TagsOnRight
\topmatter
\title{Existence of the spectral gap for elliptic operators 
}\endtitle
\author{Feng-Yu Wang}\endauthor
\affil {$_{ \text{(Department of Mathematics, Beijing Normal University,
 Beijing 100875, China)}}$\\
 $_{\text{e-mail:\ wangfy\@bnu.edu.cn}}$}\endaffil              
\def\DD{\Delta}  \def\vv{\varepsilon}
\def\<{\langle} \def\>{\rangle} \def\d{\text{d}}
 \def\nn{\nabla} \def\pp{\partial}
\def\ff{\frac} \def\ss{\sqrt} \def\si{\sigma} \def\gg{\gamma}
\def\ll{\lambda} \def\aa{\alpha} \def\bb{\beta} \def\Hess{\text{Hess}}
\def\OO{\Omega }   \def\rr{\rho} 
\thanks{Research supported in part by 
NSFC(19631060), Fok Ying-Tung Educational Foundation and Scientific Research
Foundation for Returned Overseas Chinese Scholars.
Research at MSRI is supported in part by NSF grant DMS-9701755.}\endthanks
\abstract{Let $M$ be a connected, noncompact, complete Riemannian manifold,
consider the operator $L=\DD +\nn V$ for some $V\in C^2(M)$ with 
$\exp[V]$ integrable w.r.t. the Riemannian volume element.
This paper studies the existence of the spectral gap of $L$. As
a consequence of the main result, let $\rr$ be the distance 
function from a point $o$, then the spectral
gap exists provided $\lim_{\rr\to\infty}\sup L\rr<0$ while the spectral gap
does not exist if $o$ is a pole and $\lim_{\rr\to\infty}\inf L\rr\ge 0.$ 
Moreover, the elliptic operators
on $\Bbb R^d$ are also studied.}\endabstract

\endtopmatter 
\document
\baselineskip 15pt

\head{1. Introduction }\endhead 

Let $\widetilde M$ be a $d$-dimensional, connected, noncompact, complete Riemannian 
manifold, and let $M$ be either $\widetilde M$ or an unbounded regular closed 
domain in $\widetilde M$. Next, consider $L= \DD +\nn V$ for some $V\in C^\infty(M)$  
with $Z:=\int_M\exp[V]\d x <\infty.$ Let $\d\mu=Z^{-1}\exp[V]\d x$ be defined on 
$M$. 
The spectral gap of the operator
$L$  (with Neumann boundary condition 
if $\pp M\ne \emptyset$) is characterized as

$$\ll_1 = \inf\bigg\{\ff{\mu(|\nn f|^2)}{\mu(f^2)-\mu(f)^2}:
f\in C^1(M)\cap L^2(\mu), f\ne \text{constant}\bigg\}.\tag 1.1$$
We say the spectral gap of $L$ exists if $\ll_1>0.$ From now on, we
assume that $L$ is regular in the sense that $C_0^\infty(M)$ is dense in
$W^{1,2}(M,\d\mu)$ with the Sobolev norm $\|\cdot\|_{L^2(\mu)}+\|\nn\cdot\|_{
L^2(\mu)}.$

According to  Wang$^{[12]}$ and Chen-Wang$^{[3]}$, we have $\ll_1>0$ provided
the Ricci curvature is bounded below and $\Hess_V$ is uniformly negatively
definite out of a compact domain. Actually, the recent work by the 
author [14] shows that this condition implies the logarithmic Sobolev 
inequality which is stronger than the existence of spectral gap. Moreover, 
[13] proved that the
logarithmic Sobolev inequality is equivalent to an exponential integrability
of the distance function square, which naturally refers to the negativity
of Hess$_V$ along the radial direction.

On the other hand, we know that the spectral gap
may exist if the distance function itself is exponential integrable. 
For instance, let $M=[0,\infty)$ and $L=\ff {\d^2}{\d r^2}-c\ff {\d}{\d r},
c>0$, then (see [4; Example 2.8]) $\ll_1= c^2/4>0.$ From this we may guess
that the existence of spectral gap, unlike the logarithmic Sobolev inequality,
essentially depends on the first order radial-direction derivative of $V$
rather than the second order derivative. This observation is now supported by
Corollary 1.4 in the paper.

Our study is based on the fact that $\ll_1>0$ is equivalent to $\inf\si_{ess}
(-L)>0$, where $\si_{ess}(-L)$ denotes the essential spectrum of $-L$ (with
Neumann boundary condition if $\pp M\ne\emptyset$). To see this one need only
to show that $0$ is an eigenvalue with multiplicity $1$, equivalently,
for any $f\in L^2(\mu)$ with $Lf=0$ and $\nu f|_{\pp M}=0$ ( where $\nu$ denotes
the inward unit normal vector field of $\pp M$ when it is nonempty),
one has $f$ is constant. This is a consequence of a result in Sturm [10].

Next, for fixed $o\in \widetilde M$, let $\rr(x)$ be the Riemannian distance
function from $o$. For $r>0$, let $B_r=\{x\in M: \rr(x)<r\}$ and $B_r^c=
M\setminus B_r$.
by Donnely-Li's decomposition theorem (see [5]), one has $\inf\si_{
ess}(-L)=\lim_{r\to\infty} \ll^c(r)$, where

$$ 
\ll^c(r)=\inf\{\mu(|\nn f|^2): f\in C^1(M), \mu(f^2)=1,
f=0\ \text{on}\ B_r\}.$$
Thence, $\ll_1>0$ is equivalent to $\ll^c(r)>0$ for large $r$. More precisely,
we have the following.

\proclaim{Theorem 1.1} 
$(1)$ If $\mu(B_r)>0$, then

$$\ll_1\le \ll^c(r)/\mu(B_r).\tag 1.2$$

$(2)$ Let $\ll(R)$ be the smallest positive Neumann eigenvalue
of $-L$ on $B_R$. If $\ll^c(r)>0$, then

$$\ll_1\ge \sup_{R>r} \ff {\ll^c (r)\ll (R)\mu (B_R)(R-r)^2-2\ll (R)(1-\mu(B_R))}
{2\ll(R)(R-r)^2+\ll^c(r)(R-r)^2\mu(B_R)+2\mu(B_R)}>0.\tag 1.3$$
\endproclaim

We now go to estimate the quantity $\ll^c(r)$. For $D\ge 0$, define

$$\gg(r)=\sup_{\rr(x)=r,x\notin cut(o)}L\rr(x),\ \ \ C(r)=
\int_{D+1}^r\gg(s)\d s,\ r> D.$$
Here and in what follows, the point $x$ runs over $M$.

\proclaim{Theorem 1.2} Suppose that there exists $D>0$ such that either
$\pp M\subset B_D$ or $\nu\rr \le 0$ on $\pp M\cap (B_D^c\setminus cut(o))$.
For any $r_0>D$ and positive function $f\in C[r_0,\infty)$, we have

$$\ll^c(r_0)\ge\inf_{t\ge r_0}f(t)\bigg\{\int_{r_0}^t 
\exp[-C(r)]\d r \int_r^\infty 
\exp[C(s)]f(s)\d s\bigg\}^{-1}.\tag 1.4$$
Consequently, we have $\ll_1>0$ provided there exists positive $f\in 
C[D+1,\infty)$ such that

$$\sup_{t\ge D+1}\ff 1{f(t)}\int_{D+1}^t \exp[-C(r)]\d r \int_r^\infty 
\exp[C(s)]f(s)\d s<\infty.\tag 1.5$$
\endproclaim
We remark that the assumption in Theorem 1.2 
holds if either $\pp M$ is bounded or
$o\in M$ and $M$ is convex. Especially, 
for the case $M=[0,\infty)$ and $o=0$, 
we have $C(r) =V(r)$. By [4; Theorem 2.1],
if $f'>0$ then

$$\ll_1\ge \inf_{t\ge 0}f'(t)\exp[C(t)]\bigg\{\int_t^\infty 
\exp[C(s)]f(s)\d s\bigg\}^{-1}.\tag 1.6$$
Hence, the second assertion in Theorem 1.2 can be regarded as an extension of 
[4; Theorem 2.1] to Riemannian manifolds.

\proclaim{Corollary 1.3} Under the assumption of Theorem $1.2$.
We have $\ll_1>0$ if
$$\sup_{t\ge D+1}\exp[-C(t)]\int_t^\infty \exp[C(s)]
\d s<\infty.$$ 
Consequently, $\ll_1>0$ provided
$\int_{D+1}^\infty (\gg +\vv)^+(r)\d r<\infty$ for some $\vv>0.$\endproclaim

Now, it is the time to state the result mentioned in the abstract.

\proclaim{Corollary 1.4} $(1)$ Under the assumption of Theorem $1.2$. We have 
$\ll_1>0$ provided $\lim_{r\to\infty}\sup_{\rr(x)=r, x\notin cut(o)}L\rr(x)<0$.

$(2)$ Suppose that $\pp M$ is bounded and $o$ is a pole. If $\lim_{r\to\infty}
\inf_{\rr(x)=r}L\rr(x)\ge 0$, then $\ll_1=0.$\endproclaim

\demo{Remark} 1) The first part of Corollary 1.4 follows from Corollary 1.3
directly. It was pointed out to the author by the referee that this can also be
proved by using Cheeger's inequality (c.f. [1]), the argument goes as follows.
Define the isoperimetric constant by

$$h_r= \inf\ff {A(\pp\OO\cap \text{interior}(M))}{\mu(\OO)},$$
where $\OO$ runs over all bounded open subset of $B_r^c$ and $A$ denotes the
measure on $\pp\OO$ induced by $\mu$. Then, one has $\ll^c(r)\ge h_r^2/4.$

Next, let $\nu$ be the inward normal vector field of $\pp\OO$. Noting that
$\nu\rr\le 0$ on $\pp M\cap\pp\OO$, we obtain

$$ A(\pp \OO\cap\text{interior}(M))\ge \int_{\pp\OO}\nu\rr\d A
-\int_{\pp\OO\cap
\pp M}\nu\rr\d A\ge\int_{\pp\OO}\nu\rr\d A=-\int_\OO L\rr\d\mu,$$
where $L\rr$ is understood in distribution sense in the case that cut$(o)\ne
\emptyset$. Then, under the condition we have $\lim_{r\to\infty}\ll^c(r)>0.$

2) The proof of Corollary 1.4 (2) is based on  the following upper bound
estimate (c.f. [9; Proposition 2.13]):

$$\ll_1\le \ff 1 4 \sup\{\vv^2: \mu(\exp[\vv \rr])<\infty\}.\tag 1.7$$
This estimate can be proved by taking the test function $f_n = \exp[\vv(
\rr\land n)/2]$ and then letting $n\to\infty$, refer to the proof of Theorem
3.2 below.\enddemo

The proofs of the above results are 
given in the next section, and along the same 
line, the spectral gap of elliptic operators on $\Bbb R^d$ is studied in
section 3.

\head 2. Proofs \endhead

\demo{Proof of Theorem 2.1} We prove (1.2) and (1.3) respectively.

a) The proof of (1.2) is modified from Thomas [10] which studies the upper 
bound of the spectral gap for discrete systems.

For any $\vv>0$, choose $f_\vv \in C^1(M)$ with $f_\vv|_{B_r}=0$
such that $\mu(f_\vv^2)=1$ and $\mu(|\nn f_\vv|^2)\le \vv +\ll^c(r).$
Noting that $\mu(f_\vv)=\mu(f_\vv1_{B_r^c})\le \ss{\mu(B_r^c)},$ we have 
$\mu(f_\vv^2)-\mu(f_\vv)^2\ge \mu(B_r),$ then

$$\ll_1\le \ff {\mu(|\nn f_\vv |^2)} {\mu(f_\vv^2)-\mu(f_\vv)^2} \le
\ff {\vv+\ll^c(r)} {\mu(B_r)}.$$
This proves (1.2) by letting $\vv\downarrow 0$.

b) Next, we go to prove (1.3). It suffices to show that for any $f\in C^1(M)$ 
with $\mu(f^2)=1, \mu(f)=0$ and any $R>r,$

$$\mu(|\nn f|^2)\ge \ff {\ll^c(r)\ll(R)\mu(B_R)(R-r)^2-2\ll(R)(1-\mu(B_R))}
{2\ll(R)(R-r)^2+\ll^c(r)(R-r)^2\mu(B_R)+2\mu(B_R)}.\tag 2.1$$

Let $a=\mu(f^2 1_{B_R})$. Noting that $\mu(f)=0,$
we obtain

$$\align\ff 1 {\ll(R)}\mu_{B_R}(|\nn f|^2)&\ge \mu_{B_R}(f^21_{B_R})
-\mu_{B_R}(f1_{B_R})^2= \ff a {\mu(B_R)}-\ff {\mu(f1_{B_R^c})^2} {\mu(B_R)^2}\\
&\ge \ff a {\mu(B_R)} - \ff {(1-a)(1-\mu(B_R))}{\mu(B_R)^2}
=\ff{a+\mu(B_R)-1}{\mu(B_R)^2},\endalign$$
where $\mu_{B_R}= \mu/\mu(B_R)$.
This implies

$$\mu(|\nn f|^2)\ge \ff {\ll(R)}{\mu(B_R)}\big(a+\mu(B_R)-1\big)=: g_1(a).
\tag 2.2$$

Next, define

$$h(x) = \cases 0, &\text{if}\ \rr(x)\le r,\\
1, &\text{if}\ \rr(x)\ge R,\\
\ff {\rr(x)-r}{R-r}, &\text{otherwise.}\endcases$$
Then $fh=0$ on $\{x\in M:\rr(x)=r\}.$ By the definition of $\ll^c(r)$,

$$1-a\le \mu(f^2h^2)\le \ff {\mu(|\nn(fh)|^2)}{\ll^c(r)}
\le \ff 2 {\ll^c(r)}\Big(\mu(|\nn f|^2)+\ff a {(R-r)^2}\Big).$$
Therefore,

$$\mu(|\nn f|^2)\ge \ff {\ll^c(r)} 2 -\Big(\ff {\ll^c(r)} 2
+\ff 1 {(R-r)^2}\Big)a=:g_2(a).\tag 2.3$$
By combining (2.2) with (2.3) we obtain

$$\mu(|\nn f|^2)\ge \inf_{\vv\in [0,1]}\max\{g_1(\vv), g_2(\vv)\}.\tag 2.4$$

Since $g_1(\vv)$ is increasing in $\vv$ while $g_2(\vv)$ is decreasing
in $\vv$, the above infimum is attained at

$$\vv_0= \ff {\ll^c(r)/2 +\ll(R)(1-\mu(B_R))/\mu(B_R)}
{\ll(R)/\mu(B_R)+\ll^c(r)/2+1/(R-r)^2}$$
which solves $g_1(\vv)=g_2(\vv).$ Then $\mu(|\nn f|)\ge g_1(\vv_0)=
g_2(\vv_0)$ which is equal to the right-hand side of (2.1).\qed\enddemo

\demo{Proof of Theorem 1.2} For any  $m>r_0$, let $\OO_{m}=B_m\setminus 
\bar B_{r_0}.$ Since $L$ is regular,

$$\ll^c(r_0)=\lim_{m\to\infty}\ll_0(\OO_{m}),\tag 2.5$$
where $\ll_0(\OO_{m})$ denotes the smallest eigenvalue of $-L$ on 
$\OO_{m}$ with Neumann condition on interior$(\pp M\cap\pp\OO_m)$
and Dirichlet condition on the remainder of $\pp\OO_m.$  Let $u(>0)$ be the 
corresponding eigenfunction. Define

$$F(t)=\int_{r_0}^t\exp[-C(r)]\d r\int_r^m\exp[C(s)]f(s)\d s,\ \ \ 
t\in [r_0, m].$$

We claim that there exists $c(m)>0$ such that $u(x)\le c(m)F(\rr(x))$ 
on $\OO_m$.
Actually, since $|\nn u|$ is bounded on $\OO_m$, it suffices to show that $u=0$
on $S:=\{x\in \pp \OO_m: \rr(x)=r_0\}.$ If there exists $x_0\in S$ such 
that $u(x_0)>0,$ then $x_0\in$ interior$(\pp M\cap S)$
by the boundary conditions. This
means that $\nu(x_0) =\nn \rr(x_0)$ which contradicts the 
assumption that $\nu\rr(x_0)\le 0.$

Next, let $c:=\inf_{t\ge r_0}f(t)/F(t)$, and let $x_t$ be the 
$L$-diffusion process with reflecting boundary on $\pp B_m$. By the assumption
and the It\^o's formula for $\rr(x_t)$ (see [7]), we have,
before the time $\tau:=\{t\ge 0: \rr(x_t)=r_0\}$,

$$\d\rr(x_t)=\ss 2 \d b_t +L\rr(x_t)\d t -\d L_t,\ \ \ x_0\in \OO_{m},\tag 2.6$$
where $b_t$ is an one-dimensional Brownian motion, $L\rr$ is taken to be zero
on cut$(o)$ and $L_t$ is an increasing process with support contained in $\{
t\ge 0: x_t\in \text{cut}(o)\cup \pp B_m\}.$ Noting that $L\rr(x_t)\le \gg (
\rr(x_t))$ for $x_t\notin \text{cut}(o),$ by (2.6) and It\^o's formula we obtain

$$\align \d F\circ \rr(x_t)&\le \ss 2 F'\circ \rr(x_t)\d b_t -f(x_t)\d t\\
&\le \ss 2 F'\circ \rr(x_t)\d b_t -c F\circ \rr(x_t)\d t.\endalign$$
This then implies

$$ E^x F\circ \rr(x_{t\land \tau})\le F\circ \rr(x)\exp[-ct].$$
Let $\tau'=\inf\{t\ge 0: x_t\in \overline{\pp \OO_{m}\setminus \pp M}\},$ 
we have $\tau'\le \tau$ and  $u(x_{t\land\tau'})\le
u(x_t\land \tau).$ Noting that $E^x u (x_{t\land\tau'})
=u(x)\exp[-\ll_0(\OO_{m})t],$ we obtain

$$u(x)\exp[-\ll_0(\OO_{m})t]\le c(m)E^x F\circ \rr(x_t\land\tau)\le 
c(m)F\circ\rr(x)\exp[-ct].$$
This implies $\ll_0(m)\ge c$ for any $m>r_0.$ Therefore,
$\ll^c(r_0)\ge c.\qed$\enddemo

It was pointed out by the referee that there is an equivalent analysis proof
of Theorem 1.2 (refer to [6; Lemma 1.1]). 
Let $F$ and $u$ be as in above with $\int_{\OO_m}u^2=1$, 
then $LF\circ\rr\le -c F\circ\rr $ on
$\OO_m$ in the distribution sense. Let $f = u/F$, then $f$ is bounded as was
shown in the proof of Theorem 1.2. We have $uf\nu F\le 0$ on $\pp\OO_m$ since
$\nu \rr\le 0$ on $\pp M$ and $u=0$ on $\pp\OO_m\setminus\pp M.$ Therefore,
by Green's formula, we obtain

$$\align \ll_0(\OO_m)&=-\int_{\OO_m}uLu\d\mu=-\int_{\OO_m}uL(fF)\d\mu\\
&=-\int_{\OO_m}[ufLF +uFLf+2u\<\nn f,\nn F\>]\d\mu\\
&\ge c\int_{\OO_m}[ufF+\<\nn (fF^2),\nn f\>-2u\<\nn f,\nn F\>]\d\mu
+\int_{\pp \OO_m}[uF\nu f]\d A\\
&= c+\int_{\OO_m}|\nn f|^2F\d\mu + \int_{\pp\OO_m}[u\nu u- uf\nu F]\d A\ge c.
\endalign$$ Here, we have used the fact that $u\nu u=0$ by the mixed boundary
condition.

\demo{Proof of Corollary 1.3} The prove  of the first assertion
is essentially due to [4]. Under the condition we have 

$$\int_t^\infty \exp[C(s)]\d s\le c\exp[C(t)],\ \ \ t\ge D+1$$
for some constant $c>0$. This implies (see [4; Lemma 6.1])

$$\int_t^\infty \exp[\vv s+C(s)]\d s \le \ff c {1-c\vv}\exp[C(t)+\vv t],
\ \ \ \vv\in (0,c^{-1}).$$
By taking $f(r) = \exp[r/(2c)]$ in (1.4), we prove the first assertion.

Next, if there exists $\vv>0$ such that 
$c_1:=\int_{D+1}^\infty (\gg+\vv)^+\d s <\infty.$ Let $C_\vv (r)= C(r)
-\int_{D+1}^r
(\gg +\vv)^+(s)\d s.$ Then $C_\vv'(r)=\gg(r)-(\gg +\vv)^+(r)\le -\vv.$
Therefore,

$$\align 
\exp[-C(t)]\int_t^\infty\exp[C(s)]\d s&\le \exp[-C_\vv(t)]\int_t^\infty
\exp[C_\vv(r)+c_1]\d r\\
&\le \exp[-C_\vv(t)]\int_t^\infty\exp[C_\vv(t)-\vv(r-t)+c_1]\d r\\
&=\exp[c_1]/\vv<\infty.\endalign$$
Hence $\ll_1>0$ by the first assertion.\qed\enddemo

\demo{Remark}
>From (1.3) we may derive explicit lower bounds of $\ll_1$. For instance, assume 
that $B_R$ is convex for any $R$, let $K\ge 0$ be such
that Ric$-\Hess_V\ge -K$. We have$^{[3]}$

$$\ll(R)\ge \ff {\pi^2}8K\Big\{\exp\big[KR^2/2\big]-1\Big\}^{-1}.\tag 2.7$$
Next, if $\lim_{r\to\infty}\sup \gg(r) <0$, let $\bb(r)=\inf_{s\ge r}(-\gg(s))^+,$
by taking $f(t)= \exp[\bb(r)t/2]$ in (1.4), we obtain

$$\align\ll^c(r)&\ge\inf_{t\ge r}\exp\bigg[\ff {\bb(r)t}2\bigg]
\bigg\{\int_r^t\d s\int_s^\infty\exp\bigg[\ff {\bb(r)u}2 +\int_s^u\gg(v)\d v
\bigg]\d u\bigg\}^{-1}\\
&\ge\inf_{t\ge r}\exp\bigg[\ff {\bb(r)t}2\bigg]\bigg\{\int_r^t\d s\int_s^\infty
\exp\bigg[\ff {\bb(r)u}2-\bb(r)(u-s)\bigg]\d u\bigg\}^{-1}\\
&= \ff{\bb(r)^2}4.\tag 2.8\endalign$$
Here, in the second step, we have assumed that $\bb(r)>0$ so that $\gg(v)\le
-\bb(r)$ for $v\ge r$. Then the estimate $\ll^c(r)\ge\bb(r)^2/4$ is true
for any $r$ since $\bb(r)$ is nonnegative.
The explicit lower bound of $\ll_1$ then follows from (1.3), (2.7) and (2.8).
\enddemo 

\demo{Proof of Corollary 1.4 (2)} Suppose that $\pp M\subset B_D$.
Under the polar coordinate at $o$, we have
$x=(r,\xi)$ for $r=\rr(x)$ and $\xi\in \Bbb S^{d-1},$
the $(d-1)$-dimensional unit sphere which is considered as the bundle of unit
tangent vectors at $o$. Under this coordinate, the Riemannian
volume element can be written as $\d x= g(r,\xi) \d r\d\xi$ and $\DD \rr
=\ff{\pp}{\pp r}(\log g(r,\xi))|_{r=\rr}.$ Suppose that $\lim_{r\to\infty}
\inf_{\rr(x)=r}L\rr\ge 0.$ Then, 
for any $\vv>0$ there exists $r_1>D$ such that

$$\ff{\pp}{\pp r}\log g(r,\xi)\ge -\ff \vv 2 -\ff{\pp}{\pp r}V(r,\xi),
\ \ \  r\ge r_1.$$
This implies

$$\align g(r,\xi)&\ge g(r_1,\xi)\exp\bigg[-\ff \vv 2 (r-r_1)-V(r,\xi)+
V(r_1,\xi)\bigg]\\
&\ge c\exp\bigg[-\ff \vv 2 \rr -V\bigg],\ \ \ r\ge r_1
\endalign $$
for some constant $c>0$. Therefore

$$\align \mu(\exp[\vv\rr])&\ge \int_{[D,\infty)\times \Bbb S^{d-1}}
\exp[\vv r +V(r,\xi)]g(r,\xi)\d r\d\xi\\
&\ge c\int_{[r_1,\infty)\times \Bbb S^{d-1}}
\exp\bigg[\ff \vv 2 r\bigg]\d r\d\xi=\infty.\endalign$$
By (1.7), we have $\ll_1=0.$\qed\enddemo

\demo{Remark} 1) According to the above proof, the function $\rr$ in Corollary 
1.4 (2) can be replaced by the distance 
from any bounded regular domain such that
the outward-pointing normal exponential map on the boundary induces a 
diffeomorphism. See e.g. Kumura [8] for some discussions on such manifolds.

2) In general, for any $r> D>0$, let 

$$\Xi_r^D =\{\xi\in \Bbb S^{d-1}: \exp[s\xi]|_{[0,r]}\text{\ 
is\ minimal\ and}\ \exp[s\xi]\in M, s\in [D,r]\}.$$
Then $\Xi_r^D$ is nonincreasing in $r$. Let $\nu=\d\xi$
be the standard measure on $\Bbb S^{d-1}$, the assumption of Corollary 1.4 (2)
can be replaced by: there exists $D>0$ such that

$$ \lim_{r\to\infty} \nu(\Xi_r^D)\exp[\vv r]=\infty\ \ \ \ \text{for\ any}\ 
\vv>0.\tag 2.9$$\enddemo

\head 3. Spectral gap for elliptic operators on $\Bbb R^d$\endhead

This section is a continuation of [2] and [4] in which the lower bound estimates
are studied for the spectral gap of elliptic operators on $\Bbb R^d$.

Consider the operator $L=\sum_{i,j=1}^d a_{ij}(x)\pp_i\pp_j+\sum_{i=1}^d
b_i(x)\pp_i,$ where $\pp_i=\ff\pp {\pp x_i}, a(x)$
$:=(a_{ij}(x))$ is positively 
definite, $a_{ij}\in C^2(\Bbb R^d)$ and $b_i= \sum_{j=1}^d(a_{ij}\pp_jV+\pp_ja_{
ij})$ for some $V\in C^2(\Bbb R^d)$ with $Z:= \int\exp[V]\d x<\infty.$
The specific
form of $b$ implies that $L$ is symmetric with respect to 
$\d\mu =Z^{-1}\exp[V]\d x.$ In the present setting, the spectral gap
of $L$ is described as

$$\ll_1(a,V)=\inf\{\mu(\<a\nn f,\nn f\>): f\in C^1(\Bbb R^d), 
\mu(f)=0, \mu(f^2)=1\}.\tag 4.1$$
Moreover, we assume that $L$ is regular in the sense that $C_0^\infty(
\Bbb R^d)$ is dense in $W^{1,2}(\Bbb R^d,\d\mu)$ with the Sobolev norm
$\|\cdot\|_{L^2(\mu)}+\|\ss{\<a\nn\cdot,\nn\cdot\>}\|_{L^2(\mu)}$.

Obviously, if $a\ge \aa I$ for some constant $\aa>0$,
then $\ll_1(a, V)\ge \aa\ll_1(I,V)$. From this one may transform  the present
setting to the manifold case. But this comparison only works for the case
$a$ is uniformly positively definite, and it will lead to some loss if $a$
is very different from $I$, see e.g. Examples 3.1 and 3.2 below. 
Hence, it should
be worthy to study $L$ directly as in previous sections.

Define

$$\align &\gg(r)=\sup_{|x|=r}\ff{r(\text{tr}(a(x))+\<b(x),x\>)}{\<a(x)x,x\>}
-\ff 1 r,\\
&C(r)=\int_1^r\gg(s)\d s,\ \ \ \aa(r) = \inf_{|x|=r}\ff 1 {r^2} \<a(x)x,x\>,\ \ 
r>0.\endalign$$
The main result in this section is the following.

\proclaim{Theorem 3.1} If there exists positive $f\in C[1,\infty)$ such that

$$\sup_{t\ge 1}\ff 1 {f(t)}\int_1^t\exp[-C(r)]\d r\int_r^\infty\exp[C(s)]\ff{
f(s)}{\aa(s)}\d s<\infty,\tag 3.2$$
then $\ll_1>0$.\endproclaim

\demo{Proof} For $g\in C^2(\Bbb R),$ we have

$$\align Lg(|x|)=&\ff 1 {|x|^3}\big(|x|^2\text{tr}(a(x))+|x|^2\<b(x),x\>-
\<a(x)x,x\>\big)g'(|x|)\\
&+\ff 1 {|x|^2}\<a(x)x,x\>g''(|x|),\ \ |x|>0.\endalign$$
For positive $f\in C[1,\infty)$ with $\int_1^\infty \exp[C(r)]\ff{f(r)}{\aa(r)}
\d r<\infty,$ let

$$g(t) =\int_1^t\exp[-C(r)]\d r\int_r^\infty \exp[C(s)]\ff {f(s)}{\aa(s)}\d s.$$
Then
$$ Lg(|x|)\le -f(|x|),\ \ \ |x|\ge 1.$$
Therefore, the proof of Theorem 1.2 implies that $\ll^c(1)>0$ provided
(3.2) holds, then $\ll_1>0.$\qed\enddemo

\demo{Remark} Theorem 3.1 remains true for unbounded regular domains
with bounded boundary. As for the unbounded boundary case, for the estimation
of $\ll^c(r)$, one has to consider the normal vector field induced by the 
metric $\<\pp_i,\pp_j\>=(a^{-1})_{ij},$ this will cause difficulty for general
$a$.\enddemo

For the case $M=[0,\infty),$ one has $\gg = \ff b a, \aa = a.$ 
Then, by Theorem 3.1, we have $\ll_1>0$ if there exists positive 
$f\in C[1,\infty)$ with $f'<0$ such that

$$\sup_{t>1}\ff 1 {f'(t)}\exp[-C(t)]\int_t^\infty\exp[C(s)]\ff {f(s)}{\aa(s)}
\d s<\infty.\tag 3.3$$
This is just the condition in [4; Theorem 2.1]. Therefore, Theorem 3.1 is the 
exact extension of [4; Theorem 2.1] to high dimensions.

Next, the following examples shows that Theorem 3.1 can be better than 
comparing $a$ with a constant matrix.

\demo{Example 3.1} Take $a(x)= (1+|x|^2)^\aa I, b(x)=0, \aa\ge (1+d)/2.$
Then  $L$ is regular.
It is easy to see that $V= -\aa\log (1+|x|^2)$. 
Noting that $a\ge I$, by the comparison procedure, we may consider 
the operator $\bar L=\DD -\nn V.$ But by (1.7) the spectral gap of $\bar L$ does
not exists since $\mu(\exp[\vv|x|])=\infty$ for any $\vv >0.$ Hence the comparison
procedure does not work for this example.

On the other hand, one has $\aa(r) =(1+r^2)^\aa, \gg(r) =\ff {d-1} r,
C(r)=r^{d-1}.$ Take $f=\ss t$, we obtain

$$\int_1^t\exp[-C(r)]\d r\int_r^\infty\exp[C(s)]\ff {f(s)} {\aa(s)}\d s
\le 2\int_1^t r^{3/2-2\aa}\d r\le 4f(t)$$
since $2\aa \ge 1+d$. By Theorem 3.1 we have $\ll_1>0.$\enddemo

\demo{Example 3.2} Take $a(x) = \ff 1 {|x|+1}I, V(x) =-|x|^2$ for $|x|\ge 1.$
Then the comparison procedure does not apply. Now we go to check 
the condition of Theorem 3.1. Obviously, $\aa(r)=\ff 1 {1+r}, 
\<b(x),x\>=-\ff r {(1+r)^2}-\ff {2r^2}{1+r}.$ Then

$$\gg(r)=\ff {d-1} r -\ff 1 {1+r}-2r,\ \
\exp[C(r)]=\ff {c_1 r^{d-1}}{r+1}\exp[-4r^2],
\ \ r\ge 1$$
for some $c_1>0$. Take $f(r) =r^{1-d}\exp[-r+4r^2],$ 
there exists $c_2>0$ such that

$$\ff 1 {f(t)}\int_1^t \exp[-C(r)]\d r\int_r^\infty\exp[C(s)]f(s)\d s
\le c_2 \ff 1 {f(t)}\int_1^t\exp[-r+4r^2]\ff {1+r}{r^{d-1}}\d r$$
which goes to $c_2/8$ as $t\to\infty.$ Therefore, 
Theorem 3.1 implies that $\ll_1>0$.
\enddemo

Finally, we present an upper bound estimate as (1.7).

\proclaim {Theorem 3.2} Let $\bb(r) =\sup_{|x|=r}\ff 1 {r^2}\<a(x)x,x\>,$
we have

$$\ll_1\le \ff 1 4 \sup\bigg\{\vv^2: \mu\Big(\exp\Big[\vv \int_0^{|x|}\ff 1 {
\ss{\bb(r)}}\d r\Big]\Big)<\infty\bigg\}.$$\endproclaim

\demo{Proof} Let $h(r)=\int_0^r\ff 1 {\ss{\bb(s)}}\d s.$ If $\mu(\exp[
\vv h(|x|)])=\infty,$ we go to prove $\ll_1\le \vv^2/4.$ Let $f(x)=
\exp[\ff \vv 2 (h(|x|)\land n)],\ n\ge 1.$ By (3.1), we have

$$\ll_1\le \ff {\vv^2\mu(f^2)}{4(\mu(f^2)-\mu(f)^2)}.\tag 3.4$$
Next, for any $m>1,$ choose $r_m>0$ such that $\mu(\{h(|x|)\ge r_m\})=1/m,$
we have

$$\mu(1_{\{|x|\ge r_m\}}f^2)^{1/2}\ge \ss m \mu(1_{\{|x|\ge r_m\}} f)
\ge \ss m \mu(f)-\ss m \exp[\vv h(r_m)/2].$$
Then

$$\mu(f)^2\le \Big(\ss{\mu(f^2)}\big/\ss m +\exp\big[\vv h(r_m)/2\big]\Big)^2.
\tag 3.5$$
Noting that $\mu(f^2)\to\infty$ as $n\to\infty$, by combining  (3.4) with
(3.5), we obtain

$$\ll_1\le \ff {\vv^2}{4(1-1/m)},\ \ \ m>1.$$
Therefore, $\ll_1\le \vv^2/4$ since $m$ is arbitrary.\qed\enddemo

\demo{Acknowledgement} The author is grateful to the referee whose suggestions
improved the quality of the paper. The revised version was finished during
the author's visit to MSRI in Berkeley.\enddemo

\head{References}\endhead

\ref\no 1\by Cheeger, J.\paper A lower bound for the smallest eigenvalue of the 
Laplacian\jour in ``Problems in Analysis'', Princeton Univ. Press, Princeton,
1970\endref

\ref\no 2\by Chen, M. F. and Wang, F. Y.\paper 
Estimation of the first eigenvalue of the second order elliptic operators\jour 
J. Funct. Anal. {\bf 131}(1995), 345--363\endref

\ref\no 3\by Chen, M. F. and Wang, F. Y.\paper General formula for lower
bound of the first eigenvalue on Riemannian manifolds\jour Sci. Sin. (A),
{\bf 40}(1997), 384--394\endref

\ref\no 4\by Chen, M. F. and Wang, F. Y.\paper 
Estimation of spectral gap for elliptic operators
\jour Trans. AMS. {\bf 349}(1997),
1239--1267\endref

\ref\no 5\by Donnely, H. and Li, P.\paper Pure point spectrum and negative 
curvature for noncompact manifolds\jour Duke Math. J. 
{\bf 46}(1979), 497--503\endref

\ref\no 6\by Kasue, A.\paper On a lower bound for the first eigenvalue
of the Laplace operator on a Riemannian manifold\jour
Ann. Scient. \' Ec. Norm. Sup. {\bf 17}(1984), 31--44\endref

\ref\no 7\by Kendall, W. S.\paper 
The radial part of Brownian motion on a manifold: a semimartingale property
\jour Ann. of Probab. {\bf 15}(1987), 1491--1500\endref

\ref\no  8\by Kumura, H.\paper On the essential spectrum of the Laplacian
on complete manifolds\jour J. Math. Soc. Japan. {\bf 49}(1997), 1--14\endref

\ref\no 9\by Ledoux, M.\paper Concentration of measure and logarithmic
Sobolev inequalities\jour Preprint (1997)\endref

\ref\no 10\by Strum, K. Th.\paper Analysis on local Dirichlet spaces I,
reccurrence, conservativeness, $L^p$-Liouville properties\jour
J. Riene Angrew. Math. {\bf 456}(1994), 173--196\endref

\ref\no 11\by Thomas, L. E.\paper Bound on the mass gap for finite volume 
stochastic Ising models at low temperature\jour Comm. Math. Phys. 
{\bf 126}(1989), 
1--11\endref

\ref\no 12\by  Wang, F. Y.\paper Spectral gap for diffusion processes
on noncompact manifolds\jour Chinese Sci. Bull. 
{\bf 40}(1995), 1145--1149\endref

\ref\no 13\by  Wang, F. Y.\paper Logarithmic Sobolev inequalities
on noncompact Riemannian manifolds\jour Probab. Theory Relat. Fields
{\bf 109}(1997)\endref

\ref\no 14\by Wang, F. Y.\paper Criteria of logarithmic Sobolev inequalities
for diffusion processes\jour preprint (1997)\endref
\enddocument